\documentclass[11pt]{amsart}
\usepackage{amssymb,amsmath,amsfonts,amsthm}
\usepackage[margin=1.1in]{geometry}
\usepackage{xcolor}
\usepackage{hyperref}
\hypersetup{colorlinks=true,linkcolor=blue!55!black,citecolor=blue!55!black,urlcolor=blue!55!black,
  pdftitle={A counterexample to Yau's bounded-mean-curvature embedding question},
  pdfauthor={Haoxuan Cheng},
  pdfsubject={Isometric immersions under Ricci curvature and injectivity radius bounds},
  pdfkeywords={Isometric immersions, Ricci curvature, injectivity radius, mean curvature, Weyl curvature}}
\newtheorem{theorem}{Theorem}[section]
\newtheorem{lemma}[theorem]{Lemma}
\newtheorem{proposition}[theorem]{Proposition}
\newtheorem{corollary}[theorem]{Corollary}
\theoremstyle{definition}

\theoremstyle{remark}
\newtheorem{remark}[theorem]{Remark}
\numberwithin{equation}{section}
\newcommand{\R}{\mathbb R}
\DeclareMathOperator{\Ric}{Ric}
\DeclareMathOperator{\Rm}{Rm}
\DeclareMathOperator{\Scal}{Scal}
\DeclareMathOperator{\inj}{inj}
\DeclareMathOperator{\tr}{tr}

\newcommand{\dd}{\mathrm d}
\title[Yau's bounded-mean-curvature embedding question]{A counterexample to Yau's bounded-mean-curvature\newline
embedding question}
\author{Haoxuan Cheng}
\address{School of Mathematical Sciences, Fudan University,
Shanghai 200433, China}
\email{hxcheng25@m.fudan.edu.cn}
\date{}
\subjclass[2020]{Primary 53C42; Secondary 53C20, 53C21}
\keywords{Isometric immersions, Ricci curvature, injectivity radius,
mean curvature, Weyl curvature}
\begin{document}
\begin{abstract}
Yau asked whether every complete Riemannian manifold with bounded Ricci
curvature and positive injectivity radius admits an isometric embedding into
some Euclidean space with bounded mean curvature. We construct a
counterexample: a smooth complete metric on $\R^4$ with bounded Ricci
curvature, positive injectivity
radius, and unbounded full curvature. The Gauss equation then forces every
$C^2$ isometric immersion of this metric into a finite-dimensional Euclidean
space to have unbounded mean curvature. The metric is assembled from disjoint
curvature peaks whose Ricci contractions stay uniformly bounded. A
short-geodesic curvature estimate and the Jacobi equation give the uniform
injectivity-radius bound. The construction works in every dimension at least
four, and in dimension four the Ricci curvature may be chosen to tend to zero
at infinity.
\end{abstract}
\maketitle
\section{Introduction}
\label{sec:introduction}

Nash proved that every smooth Riemannian manifold admits an isometric
embedding into some $\R^N$ \cite{Nash1956}. G\"unther later obtained a sharper
target dimension \cite{Guenther1991}. Neither theorem controls the
mean curvature of the image. In Problem~52 of his 1982 problem section, Yau
asked whether such control follows when the manifold is complete, its Ricci
curvature is bounded, and its injectivity radius has a uniform positive lower
bound \cite[p.~681]{Yau1982}. We construct a counterexample. The obstruction is
the uncontracted Gauss equation: our metric has bounded Ricci curvature and
positive injectivity radius but unbounded full curvature.

For a smooth Riemannian metric $g$, write $\Ric_g$ for its Ricci tensor,
$\Rm_g$ for its Riemann curvature tensor, and $\Scal_g$ for its scalar
curvature. Their pointwise norms are taken with respect to $g$.
We write
\[
 \inj(M,g)=\inf_{p\in M}\inj_g(p)
\]
for the global injectivity radius. Let $\delta$ denote the Euclidean metric
on $\R^4$.

\begin{theorem}\label{thm:metric}
There exists a smooth complete Riemannian metric $g$ on $\R^4$ such that
\begin{equation}\label{eq:main-properties}
 \sup_{\R^4}|\Ric_g|_g<\infty,\qquad
 \inj(\R^4,g)>0,\qquad
 \sup_{\R^4}|\Rm_g|_g=\infty.
\end{equation}
The metric can be chosen to satisfy $0.9\delta\le g\le1.1\delta$.
\end{theorem}

For an isometric immersion $F:(M^n,g)\to\R^N$, let $A$ denote its
normal-vector-valued second fundamental form and use the normalization
$H=n^{-1}\tr_g A$ for its mean curvature vector.

\begin{corollary}\label{cor:embedding}
The metric in Theorem~\ref{thm:metric} admits no $C^2$ isometric immersion
$F:(\R^4,g)\to\R^N$ with $\sup_{\R^4}|H|<\infty$, for any finite integer
$N\ge4$. In particular, it admits no such isometric embedding.
\end{corollary}

\begin{remark}\label{rem:later-yau}
Yau later asked the broader question with only a lower Ricci bound, without
assuming a positive injectivity radius \cite[\S I.H, pp.~244--245]{Yau2000};
see also \cite[\S7, p.~468]{Yau2012}. The metric in
Theorem~\ref{thm:metric} has two-sided bounded Ricci curvature, so
Corollary~\ref{cor:embedding} gives a counterexample to that version as well.
This observation does not concern variants that assume bounded full Riemann
curvature.
\end{remark}

The Michael--Simon inequalities connect Sobolev estimates on Euclidean
submanifolds with their mean curvature \cite{MichaelSimon1973}. This is one
reason that a global mean-curvature bound is a natural extrinsic condition.

Section~\ref{sec:consequences} gives several refinements.
The metric can be made arbitrarily close to $\delta$ in the global
coordinate $C^1$ norm, with an arbitrarily small Ricci bound.
It has Euclidean volume growth, uniformly over all centers and radii.
Products give examples in every dimension $n\ge4$.
A second four-dimensional construction has Ricci curvature tending to zero
at infinity, while its full curvature remains unbounded.
The scalar Gauss equation makes the implication from the theorem to the
corollary short: bounded mean curvature and a scalar curvature lower bound
force the full second fundamental form to be bounded. The uncontracted
Gauss equation then bounds $\Rm_g$. We give the argument, including the
$C^2$ regularity point, in Section~\ref{sec:gauss}.

The construction starts with a symmetric quadratic tensor $P$ on $\R^4$
whose Euclidean trace, divergence, and componentwise Laplacian vanish, and
which annihilates the radial vector. Multiplication by a radial regularized
logarithm preserves the first two identities. Its Laplacian remains bounded,
although individual second derivatives grow logarithmically. This separates
the linearized Ricci tensor from the full curvature tensor. The nonlinear
Ricci error is controlled by products of the perturbation with its second
derivatives. Placing these perturbations in disjoint balls gives the metric
in Section~\ref{sec:global}.

A uniform positive injectivity radius requires a separate argument.
The metric is uniformly close to $\delta$ in $C^1$, so every short geodesic
has a coordinate projection that increases at a uniform rate. Consequently,
it spends at most $4s$ units of arclength inside any Euclidean ball of radius
$s$. The logarithmic curvature peaks therefore have uniformly small
integrals on short geodesics. The Jacobi equation gives a bound for the
differential of the exponential map; comparison of parallel transport with
the coordinate identification then proves uniform injectivity. This is the
content of Section~\ref{sec:injectivity}.

The distinction between Ricci and full curvature bounds also appears in
the regularity theory of Riemannian metrics. In harmonic coordinates, the
Ricci tensor is an elliptic second-order operator on the metric, with
principal part $-\tfrac12 g^{ab}\partial_a\partial_b g_{ij}$
\cite[Lemma~4.1]{DeTurckKazdan1981}. Anderson proves that closed
$n$-manifolds with a uniform bound on $|\Ric|$, a uniform positive
injectivity-radius bound, and a uniform diameter bound are precompact
modulo diffeomorphisms in $C^{1,\alpha}$ for $0<\alpha<1$
\cite[Theorem~1.1]{Anderson1990}. This is a bound on metric regularity
below the second-derivative level. Our construction gives an explicit
noncompact example in which the full curvature is unbounded despite
uniform Ricci and injectivity-radius bounds.

Results under Ricci and noncollapsing assumptions naturally distinguish
low-order or integral regularity from pointwise control of the full curvature.
A Ricci lower bound and a positive injectivity-radius bound give quantitative
harmonic-coordinate control \cite{AndersonCheeger1992} and also enter
Ricci-flow smoothing without an initial full-curvature bound
\cite[Remark~1.5]{He2024}. Under two-sided Ricci bounds and noncollapsing,
Cheeger and Naber prove integral curvature estimates and codimension-four
regularity \cite{CheegerNaber2015}. Standard bounded-geometry frameworks use
the stronger hypothesis of sectional or full-curvature bounds, often together
with bounds on all covariant derivatives and a positive injectivity radius
\cite{CheegerGromovTaylor1982,Schick2001,DisconziShaoSimonett2016}.

Positive embedding results use stronger geometric control.
Petrunin proves that a complete connected smooth manifold with positive
injectivity radius and uniform bounds on every covariant derivative of its
Riemann tensor admits an isometric embedding with a uniformly thick tubular
neighborhood if and only if it has uniformly polynomial volume growth
\cite[\S2]{Petrunin2024}. Here uniform polynomial growth means that the
volume of every radius-$R$ ball is bounded by a single polynomial in $R$,
independently of its center. A uniformly thick tubular neighborhood means that
nearest-point projection is uniquely defined throughout some fixed positive
neighborhood of the image. These hypotheses include a bound on $\Rm$ itself.

Bounded mean curvature is weaker than the usual extrinsic controls. A uniform
bound on the full second fundamental form leads to
intrinsic upper-curvature and injectivity-radius estimates in suitable ambient
spaces \cite{AlexanderBishop2006} and underlies compactness theorems for
families of immersions \cite{Breuning2015}. Positive reach and $C^{1,1}$
tubular geometry form another part of this stronger theory
\cite{KapovitchLytchak2021}. Gromov obtains further scalar-curvature,
injectivity-radius, focal-radius, and mean-curvature estimates in several
geometric settings \cite{Gromov2022}.

\section{The obstruction from the Gauss equation}
\label{sec:gauss}

We use the curvature convention
\[
 R(X,Y)Z=\nabla_X\nabla_YZ-\nabla_Y\nabla_XZ-\nabla_{[X,Y]}Z,
 \qquad
 R_{ijkl}=g(R(\partial_i,\partial_j)\partial_l,\partial_k).
\]
Thus $R_{1212}$ is the sectional curvature on an orthonormal pair
$(\partial_1,\partial_2)$, and the unit sphere has positive sectional
curvature. Here $\nabla$ is the Levi-Civita connection.

\begin{proposition}\label{prop:gauss}
Let $(M^n,g)$ be a smooth Riemannian manifold whose scalar curvature is
bounded below. If $F:(M,g)\to\R^N$ is a $C^2$ isometric immersion with
bounded mean curvature vector, then $|\Rm_g|_g$ is bounded.
The curvature estimate is independent of $N$.
\end{proposition}

\begin{proof}
In a $g$-orthonormal frame, write $A_{ij}=A(e_i,e_j)$.
The Euclidean Gauss equation and its scalar contraction are
\begin{align}
 R_{ijkl}&=\langle A_{ik},A_{jl}\rangle
            -\langle A_{il},A_{jk}\rangle,\label{eq:gauss}\\
 \Scal_g&=n^2|H|^2-|A|^2.\label{eq:scalar-gauss}
\end{align}
For smooth immersions these follow from the Gauss formula
$D_X(F_*Y)=F_*(\nabla_XY)+A(X,Y)$, where $D$ is the Euclidean connection;
see \cite[Theorem~8.5]{Lee2018}.

The same identities hold for a $C^2$ immersion. To see this without
assuming third derivatives of $F$, mollify its coordinate functions on
a relatively compact coordinate ball. On the closure of a smaller ball,
$DF$ has a positive lower bound for its least singular value. The resulting
smooth maps $F_\nu$ converge in $C^2$, so for all sufficiently large $\nu$
they remain immersions there and their induced metrics are uniformly
positive definite. Write $E_\nu=DF_\nu$, so the induced metrics
$g_\nu=E_\nu^{\mathsf T}E_\nu$ converge to $g$ in $C^1$.
The orthogonal projections onto the normal spaces are
\[
 \Pi_\nu=I-E_\nu g_\nu^{-1}E_\nu^{\mathsf T},\qquad
 (A_\nu)_{ij}=\Pi_\nu\partial_i\partial_jF_\nu.
\]
Consequently $\Pi_\nu$ and the coordinate components $(A_\nu)_{ij}$
converge uniformly to their counterparts for $F$.
In the fixed coordinate basis the fully covariant curvature tensor satisfies
\[
 R_{ijkl}(g_\nu)=\tfrac12\bigl(
 \partial_i\partial_l(g_\nu)_{jk}
 +\partial_j\partial_k(g_\nu)_{il}
 -\partial_i\partial_k(g_\nu)_{jl}
 -\partial_j\partial_l(g_\nu)_{ik}\bigr)
 +Q_{ijkl}(g_\nu^{-1},Dg_\nu),
\]
where $Q$ is quadratic in the first derivatives, with coefficients
smooth in the inverse metric. Its first term converges in distributions
and $Q$ converges uniformly. Thus each $R_{ijkl}(g_\nu)$ converges to
$R_{ijkl}(g)$ in distributions. Apply the smooth Gauss equation in this
same coordinate basis, before choosing any orthonormal frame.
Its right side converges uniformly, so the limiting identity holds in
distributions and then pointwise, since both sides are continuous.
Contracting gives \eqref{eq:scalar-gauss}.

Choose constants $S,M\ge0$ such that $\Scal_g\ge-S$ and $|H|\le M$.
Equation~\eqref{eq:scalar-gauss} gives $|A|^2\le n^2M^2+S$.
By Cauchy--Schwarz in the normal space and summation over the tangent
indices, \eqref{eq:gauss} gives
\[
 |\Rm_g|_g^2\le4|A|^4.
\]
It follows that $|\Rm_g|_g\le2(n^2M^2+S)$, with no dependence on the
dimension of the normal space.
\end{proof}

\begin{proof}[Proof of Corollary~\ref{cor:embedding}, assuming Theorem~\ref{thm:metric}]
The Ricci bound gives a scalar curvature lower bound, since
$|\Scal_g|=|\tr_g\Ric_g|\le\sqrt n\,|\Ric_g|_g$. Any immersion in the
corollary would therefore have bounded intrinsic curvature by
Proposition~\ref{prop:gauss}, contradicting \eqref{eq:main-properties}.
\end{proof}

\section{A local perturbation with bounded Ricci curvature}
\label{sec:local}

All derivatives in this section are Euclidean coordinate derivatives.
For a tensor $h$, the symbols $Dh$ and $D^2h$ denote its arrays of first
and second partial derivatives, and $|h|_\delta$ is its Euclidean tensor
norm. Constants denoted by $C$ may increase from line to line; they depend
only on the fixed cutoff and the dimension, unless a dependence on a fixed
parameter is indicated.

The starting point is a nonzero algebraic curvature tensor $W$ with zero
Ricci contraction. Specify it by
\[
 W_{1212}=W_{3434}=1,\qquad
 W_{1313}=W_{2424}=-1,\qquad W_{1414}=W_{2323}=0.
\]
Extend these components by the curvature symmetries, and set
$W_{ijkl}=0$ whenever the unordered pairs $\{i,j\}$ and $\{k,l\}$
differ. The first Bianchi identity holds: components with three or four
distinct indices vanish, and the remaining cases follow from antisymmetry.
Every diagonal Ricci contraction is $1-1=0$, and every off-diagonal
contraction vanishes. Thus $W$ is an algebraic Weyl tensor.
The quadratic tensor $P_{ij}(x)=\sum_{k,l}W_{ikjl}x_kx_l$ expands as follows.

For $x=(x_1,x_2,x_3,x_4)\in\R^4$, define a symmetric covariant two-tensor
$P$ by its matrix in the standard coordinates:
\begin{equation}\label{eq:P}
 P(x)=
 \begin{pmatrix}
 x_2^2-x_3^2&-x_1x_2&x_1x_3&0\\
 -x_1x_2&x_1^2-x_4^2&0&x_2x_4\\
 x_1x_3&0&-x_1^2+x_4^2&-x_3x_4\\
 0&x_2x_4&-x_3x_4&-x_2^2+x_3^2
 \end{pmatrix}.
\end{equation}
The identities below are the reason for this choice. In these identities
the Laplacian $\Delta=\sum_{k=1}^4\partial_k^2$ acts componentwise.

\begin{lemma}\label{lem:P}
The tensor $P$ satisfies
\begin{equation}\label{eq:P-identities}
 \tr_\delta P=0,\qquad
 \sum_{i=1}^4\partial_iP_{ij}=0,\qquad
 \Delta P_{ij}=0,\qquad
 \sum_{i=1}^4x_iP_{ij}=0.
\end{equation}
Moreover, $|P|_\delta\le C|x|^2$, $|DP|_\delta\le C|x|$, and
$|D^2P|_\delta\le C$.
\end{lemma}

\begin{proof}
The diagonal entries in \eqref{eq:P} sum to zero. For the four columns,
the only nonzero divergence terms respectively sum to
$-x_1+x_1$, $-x_2+x_2$, $x_3-x_3$, and $x_4-x_4$.
The Laplacian of each diagonal entry is $2-2=0$, and every off-diagonal
entry is harmonic. The radial contractions of the four rows are
\begin{align*}
 &x_1(x_2^2-x_3^2)-x_1x_2^2+x_1x_3^2=0,\\
 &-x_1^2x_2+x_2(x_1^2-x_4^2)+x_2x_4^2=0,\\
 &x_1^2x_3+x_3(-x_1^2+x_4^2)-x_3x_4^2=0,\\
 &x_2^2x_4-x_3^2x_4+x_4(-x_2^2+x_3^2)=0.
\end{align*}
Symmetry gives the column identities in \eqref{eq:P-identities}.
The derivative bounds follow because the entries are homogeneous quadratic
polynomials.
\end{proof}

Fix a smooth function $\chi:[0,\infty)\to[0,1]$ with $\chi=1$ on
$[0,1/4]$ and $\chi=0$ on $[1/2,\infty)$.
For $0<\varepsilon\le1/8$, set
\begin{equation}\label{eq:h}
 q_\varepsilon(r)=\tfrac12\log(r^2+\varepsilon^2),\qquad
 h_\varepsilon(x)=\chi(|x|)q_\varepsilon(|x|)P(x).
\end{equation}
Each $h_\varepsilon$ is smooth, including at the origin, and vanishes
outside the closed Euclidean ball of radius $1/2$.
Near the origin the cutoff is constant. At $r=1/2$, all derivatives of
$\chi$ vanish because $\chi$ is smooth and identically zero on
$[1/2,\infty)$. Thus the radial composition is smooth across both loci.

\begin{lemma}\label{lem:bounds}
Uniformly for $0<\varepsilon\le1/8$,
\begin{align}
 \sup_{\R^4}(|h_\varepsilon|_\delta+|Dh_\varepsilon|_\delta)&\le C,
       \label{eq:C1}\\
 \sup_{\R^4}\bigl(|h_\varepsilon|_\delta|D^2h_\varepsilon|_\delta
                  +|Dh_\varepsilon|_\delta^2\bigr)&\le C,
       \label{eq:product}\\
 \tr_\delta h_\varepsilon=0,\qquad
 \sum_i\partial_i(h_\varepsilon)_{ij}=0,\qquad
 \sup_{\R^4}|\Delta h_\varepsilon|_\delta&\le C.
       \label{eq:TT}
\end{align}
Writing $r=|x|$ and $\rho=(r^2+\varepsilon^2)^{1/2}$, one also has
\begin{equation}\label{eq:second-bound}
 |D^2h_\varepsilon(x)|_\delta\le C(1+|\log\rho|)
 \quad\text{for }r\le1/4.
\end{equation}
\end{lemma}

\begin{proof}
In the ball $r\le1/4$ the cutoff is identically one, and
\[
 \partial_k q_\varepsilon=\frac{x_k}{\rho^2},\qquad
 \partial_k\partial_lq_\varepsilon
   =\frac{\delta_{kl}}{\rho^2}-\frac{2x_kx_l}{\rho^4}.
\]
The product rule and Lemma~\ref{lem:P} give
\begin{align}
 |h_\varepsilon|_\delta&\le Cr^2(1+|\log\rho|),\nonumber\\
 |Dh_\varepsilon|_\delta&\le Cr(1+|\log\rho|),\label{eq:local-bounds}\\
 |D^2h_\varepsilon|_\delta&\le C(1+|\log\rho|).\nonumber
\end{align}
Here $r\le\rho<1$, so for $r>0$ one has $|\log\rho|\le|\log r|$.
The functions $r(1+|\log r|)$ and $r^2(1+|\log r|)^2$ are bounded near
zero. At $r=0$ both $h_\varepsilon$ and $Dh_\varepsilon$ vanish.
These observations prove \eqref{eq:C1} and \eqref{eq:product} in the
inner ball. On the annulus $1/4\le r\le1/2$, all derivatives of
$q_\varepsilon$ of order at most two are uniformly bounded, so the same
claims hold there.

Let $\phi_\varepsilon(x)=\chi(|x|)q_\varepsilon(|x|)$.
It is radial, so its gradient is a scalar multiple of $x$.
The trace and divergence identities follow from
\eqref{eq:P-identities} and the product rule. In the inner ball,
homogeneity gives $\sum_kx_k\partial_kP=2P$. Consequently,
\begin{equation}\label{eq:laplacian}
 \Delta(q_\varepsilon P)
 =\left(\frac{8}{\rho^2}-\frac{2r^2}{\rho^4}\right)P
 =\frac{6r^2+8\varepsilon^2}{(r^2+\varepsilon^2)^2}P.
\end{equation}
Since $|P|_\delta\le Cr^2$, this is uniformly bounded.
The annular derivative bounds complete the proof of \eqref{eq:TT}.
\end{proof}

For a constant $a>0$, to be fixed sufficiently small, define
$g_\varepsilon=\delta+a h_\varepsilon$.

\begin{proposition}\label{prop:local}
There are constants $a_0>0$ and $C$ such that, for $0<a\le a_0$ and
$0<\varepsilon\le1/8$, the tensor $g_\varepsilon$ is a smooth metric with
\begin{equation}\label{eq:local-metric}
 0.9\delta\le g_\varepsilon\le1.1\delta,\qquad
 |\Ric_{g_\varepsilon}|_{g_\varepsilon}\le C(a+a^2).
\end{equation}
At the origin, in the standard coordinates,
\begin{equation}\label{eq:peak}
 (R_{ijkl}(g_\varepsilon)(0))=(-3a\log\varepsilon)W,
 \qquad \Ric_{g_\varepsilon}(0)=0.
\end{equation}
Uniformly for $0<a\le a_0$ and $0<\varepsilon\le1/8$, on $|x|\le1/2$,
\begin{equation}\label{eq:envelope}
 |\Rm_{g_\varepsilon}|_{g_\varepsilon}(x)
 \le C(a+a^2)\left(1+\log^+\frac1{|x|+\varepsilon}\right),
 \qquad \log^+ t=\max\{0,\log t\}.
\end{equation}
The curvature vanishes outside this ball.
\end{proposition}

\begin{proof}
Choose $a_0\le1$ so that $a_0\sup_\varepsilon\|h_\varepsilon\|_\infty
\le1/10$, using Euclidean norms. This proves the metric comparison.
The derivative of the Ricci tensor at $\delta$, applied to a symmetric
tensor $h$, is
\[
 (D\Ric)_\delta(h)_{ij}
 =\tfrac12\left(\sum_k\partial_k\partial_i h_{jk}
 +\sum_k\partial_k\partial_jh_{ik}
 -\Delta h_{ij}-\partial_i\partial_j\tr_\delta h\right).
\]
Thus \eqref{eq:TT} implies
$ (D\Ric)_\delta(h_\varepsilon)=-\tfrac12\Delta h_\varepsilon$,
which is uniformly bounded.

For completeness, the nonlinear error can be estimated directly from
\[
 \Gamma^k_{ij}=\tfrac12g^{kl}
  (\partial_i g_{jl}+\partial_j g_{il}-\partial_l g_{ij}),
 \qquad
 \Ric_{ij}=\partial_k\Gamma^k_{ij}-\partial_j\Gamma^k_{ik}
       +\Gamma^k_{kl}\Gamma^l_{ij}-\Gamma^k_{jl}\Gamma^l_{ik},
\]
where repeated indices are summed from $1$ to $4$.
Indeed,
\[
 \partial_p g^{ij}=-g^{i\alpha}(\partial_p g_{\alpha\beta})g^{\beta j}.
\]
After differentiating the Christoffel formula, the terms containing second
derivatives are of the form $g^{-1}*D^2g$, while all remaining terms are
$g^{-1}*g^{-1}*Dg*Dg$. Replacing $g^{-1}$ in the former terms by
$\delta^{-1}$ gives the linearization, and
$g^{-1}-\delta^{-1}=O(a|h_\varepsilon|_\delta)$.
The terms linear in the perturbation are exactly the preceding derivative.
In the remaining second-derivative terms the inverse metric differs from
$\delta^{-1}$ by at most $Ca|h_\varepsilon|_\delta$. Every other term is
bounded by $Ca^2|Dh_\varepsilon|_\delta^2$. Therefore
\begin{equation}\label{eq:ric-error}
 \bigl|\Ric_{g_\varepsilon}-a(D\Ric)_\delta(h_\varepsilon)\bigr|_\delta
 \le Ca^2\bigl(|h_\varepsilon|_\delta|D^2h_\varepsilon|_\delta
                 +|Dh_\varepsilon|_\delta^2\bigr).
\end{equation}
Lemma~\ref{lem:bounds} proves the Ricci bound, with the conversion to
intrinsic norms supplied by \eqref{eq:local-metric}.

At $x=0$ we have $g_\varepsilon=\delta$ and $Dg_\varepsilon=0$.
The chosen curvature convention gives
\[
 R_{1212}(0)=\partial_1\partial_2(g_\varepsilon)_{12}(0)
 -\tfrac12\bigl(\partial_1^2(g_\varepsilon)_{22}(0)
                  +\partial_2^2(g_\varepsilon)_{11}(0)\bigr).
\]
Since $P(0)=DP(0)=0$, the second derivatives of $h_\varepsilon$ there
are $(\log\varepsilon)D^2P$. The three derivatives of $P$ in this formula
are $-1,2,2$, respectively, so $R_{1212}(0)=-3a\log\varepsilon$.
For all indices, the defining formula for $P$ gives
\[
 \partial_p\partial_qP_{ab}=W_{apbq}+W_{aqbp}.
\]
Substitution into the normal-coordinate curvature formula, followed by the
curvature symmetries and the first Bianchi identity, yields
$(D\Rm)_\delta(P)=-3W$. Equivalently, putting
$b_{ij}=W_{ijij}$ for $i\ne j$, one has
$P_{ii}=\sum_{k\ne i}b_{ik}x_k^2$ and
$P_{ij}=-b_{ij}x_ix_j$, and the same component calculation gives
$R_{ijij}(0)=-3a(\log\varepsilon)b_{ij}$ while mixed unordered index pairs
give zero. This proves the full identity \eqref{eq:peak}; its Ricci
contraction is zero because that of $W$ is zero.

Finally, the coordinate curvature formula gives
\[
 |\Rm_{g_\varepsilon}|_{g_\varepsilon}
 \le C\bigl(a|D^2h_\varepsilon|_\delta
              +a^2|Dh_\varepsilon|_\delta^2\bigr).
\]
Use \eqref{eq:second-bound} in the inner ball and uniform bounds on the
annulus. Since $\rho\le r+\varepsilon\le\sqrt2\rho$, one has
$\log^+(1/\rho)\le\log\sqrt2+\log^+(1/(r+\varepsilon))$.
This proves \eqref{eq:envelope}.
\end{proof}

\section{The complete metric and curvature integrals}
\label{sec:global}

For each integer $j\ge1$, let
\[
 p_j=(4j,0,0,0),\qquad \varepsilon_j=e^{-j-3}.
\]
The closed Euclidean balls $\overline B_\delta(p_j,1/2)$ are disjoint,
their mutual distances are at least $3$, and the family is locally finite.
Define
\begin{equation}\label{eq:global-metric}
 g(x)=\delta+a\sum_{j=1}^{\infty}h_{\varepsilon_j}(x-p_j).
\end{equation}
Here $a$ is a single positive constant, independent of $j$, with
$a\le a_0$. It will be decreased once below and then fixed.

Local finiteness makes $g$ smooth. At every point at most one perturbation
is nonzero, so Proposition~\ref{prop:local} gives
\begin{equation}\label{eq:global-estimates}
 0.9\delta\le g\le1.1\delta,\qquad
 \sup_{\R^4}|\Ric_g|_g\le C(a+a^2),\qquad
 |\Rm_g|_g(p_j)=12a(j+3).
\end{equation}
Here $|W|_\delta^2=4\sum_{i<j}W_{ijij}^2=16$.
In particular, the full curvature is unbounded. If $d_g$ denotes the
Riemannian distance, metric comparison gives
\[
 \sqrt{0.9}\,|x-y|\le d_g(x,y)\le\sqrt{1.1}\,|x-y|.
\]
Thus $(\R^4,d_g)$ is complete. The Hopf--Rinow theorem implies geodesic
completeness and existence of a minimizing geodesic between any two points
\cite[Theorem~6.19 and Corollary~6.21]{Lee2018}.

The estimate \eqref{eq:C1} and the coordinate formula for the connection
give a uniform bound $|\Gamma|_\delta\le Ca$ on the Christoffel symbols.
Consequently every $g$-unit-speed geodesic $\gamma$ satisfies
\begin{equation}\label{eq:geodesic-speed}
 \frac1{\sqrt{1.1}}\le|\gamma'|_\delta\le\frac1{\sqrt{0.9}},
 \qquad |\gamma''|_\delta\le C_0a,
\end{equation}
where the derivatives use arclength and the fixed global coordinates.
Fix $a\in(0,a_0]$ small enough that $C_0a\le1/4$.
All subsequent constants and radius choices use this fixed metric.
The next two estimates hold for all short geodesics; no minimizing
assumption is made.

\begin{lemma}\label{lem:residence}
Let $\gamma:[0,L]\to\R^4$ be a $g$-unit-speed geodesic with $0<L\le1$.
For every $z\in\R^4$ and $s>0$,
\begin{equation}\label{eq:residence}
 \bigl|\{t\in[0,L]:|\gamma(t)-z|<s\}\bigr|\le\min\{L,4s\},
\end{equation}
where the outer vertical bars denote one-dimensional Lebesgue measure.
Moreover, $\gamma$ meets at most one of the balls
$\overline B_\delta(p_j,1/2)$.
\end{lemma}

\begin{proof}
Let $e=\gamma'(0)/|\gamma'(0)|_\delta$.
For $0\le t\le L$, the acceleration estimate gives
\[
 \frac{\dd}{\dd t}\langle\gamma(t),e\rangle
 =|\gamma'(0)|_\delta+
   \int_0^t\langle\gamma''(\tau),e\rangle\,\dd\tau
 \ge\frac1{\sqrt{1.1}}-C_0at\ge\frac12.
\]
Hence the inverse image of any interval of length $2s$ under
$t\mapsto\langle\gamma(t),e\rangle$ has length at most $4s$.
Since the projection of $B_\delta(z,s)$ is contained in such an interval,
this bounds the total time spent in the ball, including all possible exits
and reentries. The bound by $L$ is immediate. Finally, if $\gamma$ met two
distinct closed perturbation balls, the portion between the two intersection
points would have Euclidean length at least $3$, whereas the Euclidean length
of the entire curve is at most $L/\sqrt{0.9}\le1/\sqrt{0.9}<3$.
\end{proof}

\begin{lemma}\label{lem:integral}
There exists a constant $C_*>0$ such that every $g$-unit-speed geodesic
$\gamma:[0,L]\to\R^4$ with $0<L\le1$ satisfies
\begin{equation}\label{eq:curvature-integral}
 \int_0^L|\Rm_g|_g(\gamma(t))\,\dd t
 \le C_*L\left(1+\log\frac4L\right).
\end{equation}
\end{lemma}

\begin{proof}
If $\gamma$ meets no perturbation ball, its curvature integral is zero.
Otherwise let $p_j$ be the unique possible center and put
$\varepsilon=\varepsilon_j$. On the entire curve, including any portion
outside that ball, \eqref{eq:envelope} bounds the integrand by
\[
 C(a+a^2)\left(1+\log^+\frac1{|\gamma(t)-p_j|+\varepsilon}\right).
\]
For $u\ge0$, the superlevel set where the logarithm exceeds $u$ is
contained in $\{t:|\gamma(t)-p_j|<e^{-u}\}$.
Tonelli's theorem and Lemma~\ref{lem:residence} give
\begin{align*}
 \int_0^L\log^+\frac1{|\gamma(t)-p_j|+\varepsilon}\,\dd t
 &\le\int_0^\infty\min\{L,4e^{-u}\}\,\dd u\\
 &=L\log\frac4L+L.
\end{align*}
Absorbing the constant contribution $C(a+a^2)L$ proves the claim.
\end{proof}

\section{A uniform lower bound for the injectivity radius}
\label{sec:injectivity}

We now compare the exponential map with translation in the global
coordinates of $\R^4$. All tangent spaces are identified with $\R^4$
through these coordinates. This identification is distinct from parallel
transport; controlling their difference is part of the proof.

\begin{proposition}\label{prop:injectivity}
There exists $r>0$, independent of $p\in\R^4$, such that
\begin{equation}\label{eq:dexp-close}
 \|D\exp_p(v)-I\|_\delta<\tfrac12
 \quad\text{whenever }|v|_\delta<r.
\end{equation}
Consequently $\inj_g(p)\ge r/2$ for every $p\in\R^4$.
\end{proposition}

\begin{proof}
Fix $p$ and a nonzero vector $v\in T_p\R^4$ with
$\ell=|v|_{g_p}\le1$. Let
\[
 \gamma_v(t)=\exp_p(tv),\qquad 0\le t\le1.
\]
Its $g$-speed is $\ell$. For $w\in T_p\R^4$, vary the initial vector
through $v+sw$ and let $J(t)$ be the derivative at $s=0$ of the resulting
geodesic variation. Then
\begin{equation}\label{eq:jacobi-initial}
 J(0)=0,\qquad D_tJ(0)=w,\qquad J(1)=D\exp_p(v)w.
\end{equation}
The torsion-free property gives $D_s\partial_t\gamma=D_t\partial_s\gamma$.
Commuting these covariant derivatives and using the geodesic equation
gives the Jacobi equation
\begin{equation}\label{eq:jacobi}
 D_t^2J+R(J,\gamma_v')\gamma_v'=0;
\end{equation}
see also \cite[Theorem~10.1 and Proposition~10.10]{Lee2018}.

Let $T_v(t):T_p\R^4\to T_{\gamma_v(t)}\R^4$ be parallel transport,
and write $j(t)=T_v(t)^{-1}J(t)$ in the fixed inner product space
$(T_p\R^4,g_p)$. Define the operator
\[
 K_v(t)z=T_v(t)^{-1}
  R\bigl(T_v(t)z,\gamma_v'(t)\bigr)\gamma_v'(t).
\]
Parallel transport is an isometry, so
$\|K_v(t)\|_{g_p}\le\ell^2|\Rm_g|_g(\gamma_v(t))$.
If $\alpha(u)=\gamma_v(u/\ell)$ for $0\le u\le\ell$, then $\alpha$ is
$g$-unit-speed. Therefore, using $u=\ell t$ and
Lemma~\ref{lem:integral},
\begin{align}
 \int_0^1\|K_v(t)\|_{g_p}\,\dd t
 &\le\ell^2\int_0^1|\Rm_g|_g(\gamma_v(t))\,\dd t\nonumber\\
 &=\ell\int_0^\ell|\Rm_g|_g(\alpha(u))\,\dd u\nonumber\\
 &\le C_*\ell^2\left(1+\log\frac4\ell\right)=:\beta(\ell).
 \label{eq:beta}
\end{align}
The factor $\ell$ outside the arclength integral accounts for both the
two velocity factors in the curvature operator and the change of parameter.
Set $\beta(0)=0$; then $\beta$ is continuous at zero and increasing on
$[0,1]$.

Twice integrating \eqref{eq:jacobi} in the parallel trivialization gives
\begin{equation}\label{eq:volterra}
 j(t)=tw-\int_0^t(t-s)K_v(s)j(s)\,\dd s.
\end{equation}
Let $M=\sup_{0\le t\le1}|j(t)|_{g_p}$. When $\beta(\ell)\le1/2$,
\eqref{eq:volterra} implies
\[
 M\le|w|_{g_p}+\beta(\ell)M\le|w|_{g_p}+\tfrac12M,
 \qquad
 |j(1)-w|_{g_p}\le2\beta(\ell)|w|_{g_p}.
\]
If $Q_v$ denotes the linear map $w\mapsto j(1)$, then
\begin{equation}\label{eq:Q}
 \|Q_v-I\|_{g_p}\le2\beta(\ell).
\end{equation}

We next compare $T_v(1)$ with the coordinate identity.
In global coordinates its matrix solves
\[
 \frac{\dd}{\dd t}T_v(t)=-B_v(t)T_v(t),\qquad
 (B_v(t))^k{}_l=\Gamma^k_{il}(\gamma_v(t))(\gamma_v')^i(t),
 \qquad T_v(0)=I.
\]
Since $|\Gamma|_\delta\le Ca$ and the Euclidean length of $\gamma_v$
is at most $\ell/\sqrt{0.9}$,
\[
 \int_0^1\|B_v(t)\|_\delta\,\dd t\le Ca\ell.
\]
The integral equation for $T_v$ and Gronwall's inequality therefore give
\begin{equation}\label{eq:transport}
 \|T_v(1)-I\|_\delta\le e^{Ca\ell}-1\le C'a\ell
 \quad(\ell\le1).
\end{equation}
All constants are independent of $p$ and $v$.

By \eqref{eq:jacobi-initial}, in the fixed coordinate identifications we
have $D\exp_p(v)=T_v(1)Q_v$. Set $\kappa=\sqrt{1.1/0.9}$.
Uniform equivalence gives
$\|Q_v-I\|_\delta\le\kappa\|Q_v-I\|_{g_p}$, while the transport equation
gives $\|T_v(1)\|_\delta\le e^{Ca\ell}$. Consequently,
\begin{equation}\label{eq:dexp-estimate}
 \|D\exp_p(v)-I\|_\delta
 \le2\kappa e^{Ca\ell}\beta(\ell)+e^{Ca\ell}-1
 \le C\bigl(\beta(\ell)+a\ell\bigr).
\end{equation}
All constants here are independent of $p$, $v$, and the perturbation ball
met by the geodesic.
At $v=0$ the differential is $I$. Choose $r>0$ so small that
$\sqrt{1.1}\,r\le1$, $\beta(\sqrt{1.1}\,r)\le1/2$, and the right side
of \eqref{eq:dexp-estimate} is less than $1/2$ for
$0\le\ell\le\sqrt{1.1}\,r$. This proves \eqref{eq:dexp-close} uniformly
in $p$.

For $v,w$ in the Euclidean tangent ball $B_\delta(0,r)$, integrate the
differential along the segment from $w$ to $v$. Convexity of that ball gives
\begin{equation}\label{eq:exp-injective}
 |\exp_p(v)-\exp_p(w)-(v-w)|_\delta\le\tfrac12|v-w|_\delta.
\end{equation}
Thus $\exp_p$ is injective there. Its differential is nonsingular by
\eqref{eq:dexp-close}, so the inverse function theorem makes it a
diffeomorphism onto its image.

The $g_p$-ball $B_{g_p}(0,r/2)$ is contained in $B_\delta(0,r)$, because
$|v|_\delta\le|v|_{g_p}/\sqrt{0.9}$. To verify that this gives the usual
injectivity radius, take $v$ in the smaller ball and put $q=\exp_p(v)$.
The radial geodesic $t\mapsto\exp_p(tv)$ has constant speed $|v|_{g_p}$,
so $d_g(p,q)\le|v|_{g_p}<r/2$. By completeness there is a minimizing
geodesic from $p$ to $q$. Parametrize it affinely on $[0,1]$ and denote its
initial velocity by $u$. Then
\[
 \exp_p(u)=q,\qquad |u|_{g_p}=d_g(p,q)\le|v|_{g_p}<r/2.
\]
Both vectors lie in $B_\delta(0,r)$, so injectivity implies $u=v$.
Every radial geodesic in the smaller ball is therefore minimizing.
Conversely, any point at distance less than $r/2$ from $p$ is reached by
such a minimizing geodesic. Hence $\exp_p$ maps $B_{g_p}(0,r/2)$
diffeomorphically onto the metric ball $B_g(p,r/2)$, proving the claim.
\end{proof}

\begin{proof}[Proof of Theorem~\ref{thm:metric}]
Choose $a$ as in Section~\ref{sec:global} and take the smooth metric
\eqref{eq:global-metric}. Its metric comparison, bounded Ricci curvature,
and unbounded full curvature follow from \eqref{eq:global-estimates}; the
distance comparison immediately following that estimate proves completeness.
Proposition~\ref{prop:injectivity} supplies the uniform positive injectivity
radius.
\end{proof}

The argument also gives a sufficient condition independent of this
particular construction.

\begin{proposition}\label{prop:criterion}
Let $g$ be a smooth metric on $\R^n$ such that
$m\delta\le g\le M\delta$ for constants $0<m\le M<\infty$.
In the fixed standard coordinates, suppose that
$|\Gamma(x)(u,v)|_\delta\le G|u|_\delta|v|_\delta$ for some finite $G$.
Here $(\Gamma(x)(u,v))^k=\Gamma^k_{ij}(x)u^iv^j$ in those coordinates.
For $\ell>0$, define
\[
 I(\ell)=\sup_{\gamma}\int_0^\ell
 |\Rm_g|_g(\gamma(s))\,\dd s,
\]
where the supremum is over all $g$-unit-speed geodesics
$\gamma:[0,\ell]\to\R^n$. If $\ell I(\ell)\to0$ as
$\ell\downarrow0$, then $\inj(\R^n,g)>0$.
\end{proposition}

\begin{proof}
Metric comparison gives completeness. Put
$\beta(\ell)=\ell I(\ell)$ and
$\beta_*(s)=\sup_{0<\ell\le s}\beta(\ell)$.
For all sufficiently small $s>0$, this supremum is finite, and
$\beta_*(s)\to0$ as $s\downarrow0$. No monotonicity of $I$ or
$\beta$ is needed.

For $\ell=|v|_{g_p}>0$, use the same Jacobi variation, parallel transport
$T_v$, and map $Q_v$ as above. Reparametrization by arclength gives
$\int_0^1\|K_v(t)\|_{g_p}\,\dd t\le\beta(\ell)$.
The Volterra equation \eqref{eq:volterra} therefore yields
$\|Q_v-I\|_{g_p}\le2\beta(\ell)$ whenever $\beta(\ell)\le1/2$.
The coordinate length of the geodesic is at most $\ell/\sqrt m$, so
\[
 \|T_v(1)\|_\delta\le e^{G\ell/\sqrt m},\qquad
 \|T_v(1)-I\|_\delta\le e^{G\ell/\sqrt m}-1.
\]
With $\kappa=\sqrt{M/m}$, the factorization $D\exp_p(v)=T_v(1)Q_v$
now gives
\[
 \|D\exp_p(v)-I\|_\delta
 \le2\kappa e^{G\ell/\sqrt m}\beta(\ell)
       +e^{G\ell/\sqrt m}-1.
\]
Choose $s>0$ such that $\beta_*(s)\le1/2$ and the right side, with
$\beta(\ell)$ replaced by $\beta_*(s)$ and $\ell$ by $s$, is less
than $1/2$. Set $r=s/\sqrt M$ and $\rho=\sqrt m\,r/2$.
Then the differential is uniformly within $1/2$ of $I$ on every
$B_\delta(0,r)\subset T_p\R^n$, including $v=0$.
Integration along straight segments in that convex ball proves
injectivity exactly as in \eqref{eq:exp-injective}; the differential
is also nonsingular there. Since
$B_{g_p}(0,\rho)\subset B_\delta(0,r)$, completeness supplies minimizing
geodesics whose initial velocities stay in the same coordinate ball.
Injectivity forces them to coincide with the radial geodesics, as in
the preceding proof. Hence $\inj_g(p)\ge\rho$ uniformly in $p$.
\end{proof}

This criterion retains both the global coordinate hypothesis and the
connection bound. Small curvature integrals alone are not the hypothesis.
The logarithmic peaks above satisfy the stronger estimate
$I(\ell)\le C_*\ell(1+\log(4/\ell))$ for $0<\ell\le1$.

\section{Further properties of the examples}
\label{sec:consequences}

\begin{corollary}\label{cor:quantitative}
For the metric \eqref{eq:global-metric}, every $C^2$ isometric immersion
$F:(\R^4,g)\to\R^N$ satisfies
\[
 |H(p_j)|\ge\sqrt{\frac{3a}{8}}\sqrt{j+3}\qquad(j\ge1).
\]
The points $p_j$ and the constant in this bound do not depend on $F$ or $N$.
\end{corollary}

\begin{proof}
Equation~\eqref{eq:peak} gives
\[
 \Ric_g(p_j)=0,\qquad \Scal_g(p_j)=0,\qquad
 |\Rm_g|_g(p_j)=12a(j+3).
\]
At these points the scalar Gauss equation is $|A|^2=16|H|^2$.
The estimate $|\Rm_g|_g\le2|A|^2$ from Proposition~\ref{prop:gauss}
therefore gives $12a(j+3)\le32|H(p_j)|^2$.
\end{proof}

\begin{corollary}\label{cor:dimensions}
For every integer $n\ge4$ and every $\eta>0$, there exists a smooth
complete metric $g_n$ on $\R^n$ such that
\[
 \|g_n-\delta\|_{C^1(\delta)}<\eta,\qquad
 \sup_{\R^n}|\Ric_{g_n}|_{g_n}<\eta,\qquad
 \inj(\R^n,g_n)>0,\qquad
 \sup_{\R^n}|\Rm_{g_n}|_{g_n}=\infty.
\]
Here $\|h\|_{C^1(\delta)}=\sup_x(|h(x)|_\delta+|Dh(x)|_\delta)$
uses the fixed standard coordinates. There is no $C^2$ isometric immersion
of this metric into any finite-dimensional Euclidean space with bounded
mean curvature. Moreover, for constants $c_n,C_n>0$ depending only on $n$,
\[
 c_n R^n\le\operatorname{Vol}_{g_n}B_{g_n}(p,R)\le C_n R^n
 \qquad(p\in\R^n,\ R>0).
\]
\end{corollary}

\begin{proof}
In dimension four, \eqref{eq:C1} and disjointness of the supports give
$\|g-\delta\|_{C^1(\delta)}\le Ca$, and
\eqref{eq:global-estimates} gives $\sup|\Ric_g|_g\le C(a+a^2)$.
Choose the positive parameter $a$ small enough to make both bounds less
than $\eta$, in addition to the choice in Section~\ref{sec:global}.
The curvature peaks remain unbounded for every such fixed $a>0$,
and Proposition~\ref{prop:injectivity} applies.

For $n>4$, take $g_n=g\oplus\delta$ on
$\R^4\times\R^{n-4}$. The connection and curvature split into the two
factors. Thus the Ricci and full curvature norms equal those of the
four-dimensional factor, and the coordinate $C^1$ bound is unchanged.
The exponential map is
\[
 \exp^{g_n}_{(p,y)}(v,w)=(\exp^g_p(v),y+w).
\]
For $| (v,w)|_{g_n}<i_0$, where $0<i_0<\inj(\R^4,g)$,
its differential is block diagonal and nonsingular. If two such vectors
have the same image, their Euclidean components agree, and injectivity of
$\exp^g_p$ then makes their four-dimensional components agree as well.
Thus the product exponential map is injective there. Product completeness
and the minimizing-geodesic argument used in
Proposition~\ref{prop:injectivity} give $\inj(\R^n,g_n)\ge i_0$.
Apply Proposition~\ref{prop:gauss} in dimension $n$ to exclude bounded
mean curvature. This uses the $n$-dimensional immersion itself.

Finally, put $m=0.9$ and $M=1.1$. The comparison
$m\delta\le g_n\le M\delta$ gives
\[
 B_\delta(p,R/\sqrt M)\subset B_{g_n}(p,R)
 \subset B_\delta(p,R/\sqrt m),\qquad
 m^{n/2}\le\sqrt{\det g_n}\le M^{n/2}.
\]
Integration proves the volume estimates, with
$c_n=\omega_n(m/M)^{n/2}$ and $C_n=\omega_n(M/m)^{n/2}$,
where $\omega_n$ is the volume of the Euclidean unit ball.
\end{proof}

The dimension threshold comes from the algebra of curvature.
In dimension two, $\Ric=K g$. In dimension three, the curvature tensor is
determined by its Ricci contraction through
\[
 R_{ijkl}=\Ric_{ik}g_{jl}+\Ric_{jl}g_{ik}
 -\Ric_{il}g_{jk}-\Ric_{jk}g_{il}
 -\tfrac12\Scal(g_{ik}g_{jl}-g_{il}g_{jk}).
\]
Thus, in an orthonormal basis diagonalizing $\Ric$, the sectional
curvatures are
\[
 K_{ij}=\Ric_{ii}+\Ric_{jj}-\tfrac12\Scal,
\]
and the displayed tensor identity makes the off-diagonal entries of the
curvature operator vanish. Thus a Ricci norm bound implies a full curvature
bound in these dimensions. This explains why the present construction starts
in dimension four. Dimensions two and three require a different mechanism.

\begin{proposition}\label{prop:decay}
For every $\eta>0$, there is a smooth complete metric $\widetilde g$ on
$\R^4$ satisfying the conclusions of Corollary~\ref{cor:dimensions}
for $n=4$, and in addition
\[
 |\widetilde g(x)-\delta|_\delta+|D\widetilde g(x)|_\delta
 +|\Ric_{\widetilde g}(x)|_{\widetilde g}\longrightarrow0
 \qquad\text{as }|x|\longrightarrow\infty.
\]
Its full curvature is unbounded. Every $C^2$ isometric immersion into a
finite-dimensional Euclidean space satisfies
\[
 |H(p_j)|\ge\sqrt{\frac{3a}{8}}(j+1),
\]
where $a>0$ is the fixed parameter chosen in the construction below.
\end{proposition}

\begin{proof}
Keep the same centers $p_j$, but set
\[
 a_j=\frac{a}{j+1},\qquad
 \widetilde\varepsilon_j=e^{-(j+1)^3},\qquad
 \widetilde g=\delta+\sum_{j\ge1}
 a_jh_{\widetilde\varepsilon_j}(x-p_j).
\]
Here $j\ge1$, so $\widetilde\varepsilon_j\le1/8$ and $0<a_j\le a$.
Choose $a>0$ small enough for the preceding estimates and the two
bounds by $\eta$. Local finiteness and metric comparison again prove
smoothness and completeness. On the $j$th perturbation ball,
\[
 |\widetilde g-\delta|_\delta+|D\widetilde g|_\delta\le Ca_j,
 \qquad |\Ric_{\widetilde g}|_{\widetilde g}\le C(a_j+a_j^2).
\]
Both expressions vanish outside the balls. Any sequence tending to
infinity that meets these balls must have ball indices tending to
infinity, proving the asserted decay.
At their centers, \eqref{eq:peak} gives
\[
 \Ric_{\widetilde g}(p_j)=0,\qquad
 R_{1212}(\widetilde g)(p_j)=3a(j+1)^2,\qquad
 |\Rm_{\widetilde g}|_{\widetilde g}(p_j)=12a(j+1)^2.
\]
This proves curvature divergence and, by the proof of
Corollary~\ref{cor:quantitative}, the mean curvature lower bound.

It remains to check the injectivity radius. Since $a_j\le a$,
the Christoffel and acceleration bounds in
\eqref{eq:geodesic-speed} hold with the same constants.
The residence-time argument therefore still applies to every geodesic
of length at most one. The amplitude dependence made explicit in
\eqref{eq:envelope} gives on each ball
\[
 |\Rm_{\widetilde g}|_{\widetilde g}(x)
 \le C(a+a^2)\left(1+\log^+
       \frac1{|x-p_j|+\widetilde\varepsilon_j}\right).
\]
The proof of Lemma~\ref{lem:integral} consequently gives the uniform estimate
\[
 \widetilde I(\ell):=\sup_\gamma\int_0^\ell
 |\Rm_{\widetilde g}|_{\widetilde g}(\gamma(s))\,\dd s
 \le C_*\ell\left(1+\log\frac4\ell\right),
 \qquad 0<\ell\le1,
\]
where the supremum is over unit-speed geodesics of length $\ell$.
Hence $\ell\widetilde I(\ell)\to0$. The global coordinate comparison
and Christoffel bound are already uniform, so
Proposition~\ref{prop:criterion} gives $\inj(\R^4,\widetilde g)>0$.
The volume bounds follow from metric comparison as above.
\end{proof}

Proposition~\ref{prop:decay} gives decay of the metric, its first coordinate
derivatives, and its Ricci tensor. Its full curvature remains unbounded and
diverges along the centers $p_j$, so this is not asymptotic flatness in a
curvature-decay sense.

\section*{Acknowledgements}
The author thanks Professor Bobo Hua for his support.

\section*{AI usage}
The author used OpenAI's Codex in this work. This assistance included proposing and refining candidate metric constructions and proofs, performing computational checks, and reviewing arguments through separate AI agents. The author guided the research direction, checked the mathematical statements and proofs, and carried out the final revisions. The author takes full responsibility for the content of the paper.

\bibliographystyle{amsplain}
\bibliography{references}
\end{document}